\documentclass[12pt]{article}
\usepackage{amsmath, amsfonts, amsthm}

\newtheorem{rem}{Remark}

\newcommand{\bg}{\overline{\gamma}}
\newcommand{\tg}{\widetilde{\gamma}}

\newcommand{\mN}{\mathcal{N}}
\newcommand{\mI}{\mathcal{I}}
\newcommand{\bh}{\overline{h}}
\newcommand{\bp}{\overline{\pi}}

\newcommand{\Rg}       {{\hbox{I\kern-.22em\hbox{R}}}}
\newcommand{\Pg}       {{\hbox{I\kern-.22em\hbox{P}}}}
\newcommand{\Eg}       {{\hbox{I\kern-.22em\hbox{E}}}}

\newcommand{\YY}{{\cal Y}}

\newtheorem{theorem}{Theorem}

\newtheorem{proposition}[theorem]{Proposition}
\newtheorem{corollary}[theorem]{Corollary}

\theoremstyle{definition}

\newcommand{\argmin}{\mathop{\mathrm{argmin}}}
\newcommand{\cov}{\mathop{\mathrm{cov}}\nolimits}
\newcommand{\var}{\mathop{\mathrm{var}}\nolimits}
\newcommand{\E}{\mathop{\mathbb{E}}\nolimits}

\title{Filtering problems with exponential criteria for general Gaussian signals}
\author{
M.L.~Kleptsyna \footnote{Author for correspondence.
Fax: +33 2 43 83 35 79.}\\
\vspace{-2mm} \small\sl
Laboratoire de Statistique et Processus / Universit\'e du Maine\\
\vspace{-2mm}
\small\sl Av. Olivier Messiaen, 72085 Le Mans, Cedex 9, France \\
\small\sl e-mail: Marina.Kleptsyna@univ-lemans.fr\\
\mbox{\hspace{100mm}}\vspace{-5mm} \and \protect{\newline}
\vspace{-1mm}
A.~Le Breton \\
\vspace{-2mm} \small\sl
Laboratoire Jean Kuntzmann / Universit\'e J. Fourier\\
\vspace{-2mm}
\small\sl BP 53, 38041 Grenoble Cedex 9, France  \\
\small\sl e-mail: Alain.Le-Breton@imag.fr\\
\mbox{\hspace{100mm}}\vspace{-5mm} \and \protect{\newline}
\vspace{-1mm} M.~Viot\\
\vspace{-2mm} \small\sl
Laboratoire Jean Kuntzmann / Universit\'e J. Fourier\\
\vspace{-2mm} \small\sl BP 53, 38041 Grenoble Cedex 9, France}

\date{}

\begin{document}
\maketitle
\begin{abstract} The explicit solution of the discrete time filtering problems with exponential criteria for a general
Gaussian signal is obtained through an  approach  based on a conditional Cameron-Martin type formula.
This key formula
is derived for  conditional expectations
of exponentials of some quadratic forms of Gaussian sequences. The formula
involves  conditional expectations
and  conditional covariances  in some auxiliary optimal risk-neutral filtering problem
which  is used in the proof. Closed form recursions of  Volterra
type for these ingredients are provided. Particular cases for which the results can be further
elaborated are investigated.
\end{abstract}

\vspace*{4mm}


\vspace*{0.25cm}


\section{Introduction}
The linear exponential Gaussian (LEG for short) filtering problem,
\textit{i.e.}, with an exponential cost criteria (see the definition
(\ref{LEGdef}) below), and  the so called risk-sensitive (RS for short) filtering
problem (see \cite{collings} and  the statement (\ref{rde})
below) have been given a great deal of interest over the last
decades. Numerous results have been already reported in specific
models, specially around  Markov models, but, as far as we know,
without exhibiting  the relationship between these two problems.
See, \textit{e.g.}, Whittle \cite{whittle4}-\cite{whittle3}, Speyer  \textit{et
al.} \cite{speyer}, Elliott \textit{et al}. \cite{elliott4}, \cite{dey1},
\cite{elliott1}  and  Bensoussan and van Schuppen
\cite{bensoussan} for contributions on this subject and related LEG and RS control problems. Therein
the notion of  ``information state" has been introduced
without any clear
 probabilistic meaning for auxiliary  processes which are involved,
 even in the Gauss-Markov case.
Moreover, the method proposed in \cite{elliott4} does not work in
a non Markovian situation.
 In our paper \cite {AMM}, we have  solved the LEG and RS filtering problems  for general
 Gaussian signal processes in continuous time  and in the particular setting where the  functional in the exponential is a \textit{singular} quadratic functional. Moreover we have  proved that actually in this case the solutions coincide.
 In our  paper \cite{AMM2} we have  solved the LEG and RS filtering problems for Gauss-Markov processes but with a \textit{nonsingular} quadratic functional in the exponential.
 In this setting   we have proposed  an example to show that the solutions  may be different.
On the other hand,
  the  general solution for the
optimal  risk-neutral linear filtering problem and
 a Cameron-Martin type formula for  general Gaussian sequences have been obtained in \cite{mkalbmcr3}. It seems natural
to use the approach proposed in \cite{mkalbmcr3} and \cite{AMM}  to  derive the solution of the
LEG and RS filtering problems for  general Gaussian signals in discrete time setting, to precise their link and also to give a
probabilistic  interpretation for the ingredients of the ``information state".

In the present paper we are interested in the explicit solution of the Linear Exponential Gaussian (LEG)  and Risk Sensitive (RS)  filtering problems for general Gaussian signals. Namely we deal with a signal-observation model
$
(X_t,Y_t)_{t\ge 1},
$
where the signal $X=(X_t)_{t\ge 1}$ is an arbitrary Gaussian sequence with mean
$m=(m_t, t\geq 1)$ and covariance   $K =(K(t,s), t\geq 1,
s\geq 1)$, \textit{i.e.},
$$
\Eg X_t=m_t,\quad\Eg (X_t-m_t)(X_s-m_s)=K(t,s)\,, \quad
t\geq 1\,,\; s\geq 1\,,
$$
and, for some sequence $A=(A_t,\, t\ge 1)$ of the real numbers, the observation process $Y=(Y_t,\, t\ge 1)$ is given by
\begin{equation}\label{observ}
Y_t= A_t X_t +\varepsilon_t,
\end{equation}
where $\varepsilon=(\varepsilon_{t})_{t\ge 1}$ is a sequence of i.i.d. $\mN(0,1)$ random variables and $\varepsilon$ and $X$ are independent.

Suppose that only $Y$ is observed and for a given real number $\mu$ and a fixed sequence $(Q_{t})_{t\ge 1}$ of nonnegative real numbers, one wishes
to minimize with respect to $h:\,h_t\in \YY_{t}, t\ge 1$  the quantity:
\begin{equation}\label{LEGcrit}
 \E \mu \exp \left\{\frac{\mu}{2} \sum_{t=1}^T (X_t-h_t)^2 Q_t \right\},
\end{equation}
where $({\YY}_t)$ is the natural filtration of $Y$, \textit{i.e.},
${\YY}_t=\sigma(\{Y_u\, ,\, 1\leq u\leq t\})$ and $h_t\in{\YY}_t$ means that  $h_t$ is ${\YY}_t$-measurable.

 Note
that, according to the sign of the real parameter $\mu$, there are two different
cases for this linear exponential Gaussian (LEG)
filtering problem (the terminology is taken from the linear exponential Gaussian optimal control problem) :
\begin{itemize}
  \item$\mu < 0 ,$ called \textit{risk-preferring} filtering problem,

  \item  $\mu > 0, $ called the \textit{risk-averse} filtering problem.

\end{itemize}

It is well known (see, \textit{e.g.}, \cite{speyer} for the
 Markov case)  that the solution to this problem is
not the conditional expectation of $X_t$ given the $\sigma$-field
${\cal Y}_t$. Our first aim is to show that the solution
can be completely explicited~: the characteristics of the optimal
solution are obtained as the solution of a closed form system of
Volterra type equations which actually  reduce to the
 equations known also for the RS setting  when the signal process $X$ is
Gauss-Markov (see, \textit{e.g.}, \cite{elliott2}). Our second aim is to give
the probabilistic interpretation of this optimal solution in terms of an
auxiliary  \textit{risk-neutral } filtering problem. Actually, we extend the
filtering approach initiated in \cite{mkalbmcr3} and \cite{AMM} for
one-dimensional processes, to obtain a conditional Cameron-Martin
type formula for the {\em conditional Laplace transform} of a
quadratic functional of the involved process. Namely, we give an
explicit representation for the random variable
\begin{equation}\label{def:It}
\mI_T= \E \left(\left.\exp\left\{\frac{\mu}{2} \sum_{s=1}^T (X_s-h_s)^2 Q_s \right\} \right/ {\YY}_{T}\right),
\end{equation}
where $h_s\in \YY_{s}, \, s\ge 1$.

\medskip The paper is organized as follows. In Section \ref{LEG} we derive
the solution of the LEG filtering problem : explicit recursive
 equations, involving the covariance function of
the filtered process, are obtained. In particular,  in Section \ref{CMC}, an appropriate auxiliary risk-neutral filtering problem is
matched to that of deriving  the key
 Cameron-Martin type formula. The solution of this auxiliary  filtering problem is discussed
in Section \ref{AFP}. In Section \ref{PC}  we investigate some specific cases
where the results can be further elaborated. In Section \ref{disc} we discuss  the relationship between LEG and
RS filtering problems. Section \ref{interpret}
is devoted to the interpretation for the ingredients of the ``information state''.
Finally, Sections \ref{complements} and \ref{againpart} are devoted to a more general case, namely when the particular structure of the observation sequences $(Y_t)_{t\ge 1}$ is not specified.

\section{Solution of the LEG filtering problem}\label{LEG}
Let us introduce the following condition $(C_{\mu})$:
\begin{enumerate}
  \item[$(C_{\mu})$]
the  equation
\begin{equation}\label{GAMMABAR}
\bg(t,s)=K(t,s)-\sum_{l=1}^{s-1} \bg(t,l)\bg(s,l) \, \frac{S_l}{1+S_l \bg_l},\quad S_l=A_l^2-\mu Q_l
\end{equation}
has a unique and bounded solution on $\{(t,s):1\le s \le t \le T\}$,   such that $\bg_l=\bg(l,l) \ge 0,\, l\ge 1$  and moreover
$$ \displaystyle{1 +
 S_{l}\bg_l}> 0,\, l\ge 1.
 $$
\end{enumerate}
\begin{rem}

Notice that for all  $\mu$ \textbf{negative} the condition $(C_{\mu})$ is satisfied  and if $\mu $ is  \textbf{positive}, the condition $(C_{\mu})$ is satisfied for $\mu$  sufficiently small, for example, those such  that for any
  $ t \le T \, A_{t}^{2}-\mu Q_{t}$ is nonnegative.
\end{rem}
The first result is the following
\begin{theorem}\label{LEGsol}
Suppose that the condition $(C_{\mu})$ is satisfied. Let $(\bh_t)_{t\ge 1}$ be the solution of the following equation:
\begin{equation}\label{hbar}
\bh_t=m_t + \sum_{l=1}^t A_l\bg(t,l) (Y_l-A_l \bh_l),
\end{equation}
where  $\bg=(\bg(t,s), 1\le s \le t \le T)$ is the unique solution of equation \eqref{GAMMABAR}.\\
Then $(\bh_t)_{t\ge 1}$ is the solution of the LEG filtering problem,
 \textit{i.e.},
\begin{equation}\label{LEGdef}
\bh=\argmin_{h:\,h_t\in \YY_{t}, t\ge 1} \E \mu \exp \left\{\frac{\mu}{2} \sum_{t=1}^T (X_t-h_t)^2 Q_t \right\}.
\end{equation}
Moreover, the corresponding optimal risk is given by
$$
\E \mu \exp \left\{\frac{\mu}{2} \sum_{t=1}^T (X_t-\bh_t)^2 Q_t \right\}
=\mu \prod_{t=1}^T \left[\frac{1+S_t \bg_t}{1+A_t^{2} \bg_t}\right]^{-1/2}.
$$

\end{theorem}

Theorem \ref{LEGsol} is a direct consequence of
results of Section \ref{CMC}. Its proof will be given at the end of
Section \ref{CMC}.

\begin{rem}

\begin{itemize}
  \item Note that equation \eqref{hbar} is really recursive equation and it  can be rewritten in the equivalent form:
 $$
 \bh_t=\frac{1}{1+A_t^2\bg_{t}}\left[m_t + \sum_{l=1}^{t-1} A_l\bg(t,l) (Y_l-A_l \bh_l)+ A_t \bg_{t}Y_{t}\right],
$$

  \item It is worth  emphasizing that taking  $\mu=0$ in  equation
(\ref{GAMMABAR}), one gets through equation (\ref{hbar}) the solution $\bar{h}$ of  the
\textit{risk-neutral} filtering problem  of the signal $X$ given
the observation  $ Y$, \textit{i.e.}, $\bar{h}_{t}= \Eg ( X_{t} / {\cal{ Y }}_{t})$  (see, \textit{e.g.}, \cite{mkalbmcr3}).
\end{itemize}
\end{rem}

\subsection{Conditional version of a Cameron-Martin formula}\label{CMC}
The proof of Theorem~\ref{LEGsol} is based on the conditional version of the Cameron--Martin formula which provides
the conditional expectation $\mI_t$ defined by \eqref{def:It}. Let
\begin{equation}\label{def:Jt}
J_t=\exp\left\{ -\frac{1}{2}\sum\limits_{s=1}^t (X_s- h_s)^2 Q_s \right\}.
\end{equation}
Then $\mI_t=\pi_t (J_t)$,
where for any random variable $\eta$ such that
$\Eg |\eta|<+\infty$, the notation
$\pi_t(\eta)$ is used for the conditional expectation of $\eta$ given the
$\sigma$-field
${\cal
Y}_t=\sigma(\{Y_s\, ,\, 1\leq s\leq t\}),$
$$
\pi_t(\eta) = \Eg(\eta/{\cal Y}_t)\,.
$$

\begin{proposition}\label{p:CM}
Suppose that the condition $(C_{\mu})$ is satisfied. Let $(\bg(t,s),\, 1\le s \le t \le T)$ be the solution of equation \eqref{GAMMABAR} and  $(Z_t^h,\,t\ge 1)$ be the solution of the following equation
\begin{equation}\label{eqzh}
Z_t^h = m_t -\sum_{l=1}^{t-1}\bg(t,l) \frac{\mu Q_l}{1+S_l \bg_l} (h_l - Z_{l}^h) + \sum_{l=1}^{t-1}\bg(t,l) \frac{A_l}{1+S_l \bg_l} (Y_l - A_l Z_{l}^{h}).
\end{equation}
Then the following representation of the random variable $\mI_T$ defined by \eqref{def:It}  holds for any $T\ge 1$:
$$
\mI_T=\prod_{t=1}^T \left[\frac{1+S_t \bg_t}{1+A_t^{2} \bg_t}\right]^{-1/2} \times \exp\left\{\frac{\mu}{2} Q_t \frac{1+A_t^2 \bg_t}{1+ S_t \bg_t} \times \left[ h_t - \frac{Z_{t}^h+A_t \bg_t Y_t}{1+ A_t^2 \bg_t } \right]^{2}\right\} \times \mathcal{M}_T,
$$
where   $(\mathcal{M}_T)_{T\ge 1}$ is a martingale defined by :
\begin{multline}\label{martin}
\mathcal{M}_T = \prod_{t=1}^T \left[\frac{(1+A_t^2 \gamma_t)}{1+A_t^{2} \bg_t}\right]^{1/2} \exp\left\{ \frac{A_t}{1+A_t^2 \bg_t} \, (Z_{t}^h - \pi_{t-1}(X_t) ) \nu_t - \right.
\\
\left. -\frac{1}{2} \cdot \frac{A_t^2}{1+A_t^2 \bg_t} \, (Z_{t}^h - \pi_{t-1}(X_t))^2
- \frac{1}{2} \cdot \frac{A_t^2 (\gamma_t - \bg_t)\cdot \nu_t^2}{(1+A_t^2 \bg_t) (1+A_t^2 \gamma_t)}
\right\},
\end{multline}
in terms of the innovation sequence $(\nu_t)_{t \ge 1}$:
$$
\nu_t=Y_t - A_{t}\pi_{t-1}(X_t);\quad \pi_{t-1}(X_t)=\E (X_t / \YY_{t-1}),
$$
and of the variances of one-step prediction errors $(\gamma_t)_{t \ge 1}$:
$$
\quad \gamma_t = \E (X_t-\pi_{t-1}(X_t))^2.
$$
\end{proposition}
\begin{rem}
\begin{enumerate}
  \item
  The probabilistic interpretation  of the auxiliary processes $(Z_t^h)$ and $(\bg_t)_{t\ge 1}$ appearing in the Proposition~\ref{p:CM} will be clarified below.

  \item Proposition~\ref{p:CM} reduces to the ordinary Cameron-Martin type formula (\textit{cf.} Theorem~1 \cite{mkalbmcr3}) for $h\equiv 0$ when $A_t=0,\,l\ge 1$ and hence $X$ and $Y$ are independent.
\end{enumerate}
\end{rem}
\paragraph{Proof of Proposition~\ref{p:CM}}
We will prove Proposition~\ref{p:CM} for $\mu<0$, namely $\mu=-1$. Then we can replace $Q$ by $-\mu Q$ and the statement of Proposition~\ref{p:CM} is still valid because of the analytical properties of the involved functions.

The proof of  Proposition~\ref{p:CM} for $\mu=-1$ will be separated into two steps.

\textbf{I.} (Actually it is the discrete time analog for the general filtering theorem.) Since $h_t\in \YY_{t},\, t\ge 1,$ in the proof  we can suppose that $h$ is a deterministic function.
First of all, we claim that for $J_t$, defined by \eqref{def:Jt}
\begin{equation}\label{pitJ}
\pi_t (J_t) = \left. \frac{\pi_{t-1}(J_t \beta_t^y)}{\pi_{t-1}(\beta_t^y)} \right|_{y=Y_t},
\end{equation}
where $\beta_t^y = \exp (A_t X_t y - \frac{1}{2} A_t X_t^2)$.

Indeed, let us introduce the new probability measure $\hat{\Pg}$, defined by
$$
\frac{d\hat{\Pg}}{d\Pg} = \exp (-A_t X_t \varepsilon_t - \frac{1}{2} A_t^2 X_t^2).
$$

The classical Bayes formula gives that
$$
\pi_t(J_t) = \frac{\hat{\pi}_t (J_t \exp (A_t X_t \varepsilon_t + \frac{1}{2} A_t^2 X_t^2))}{\hat{\pi}_t( \exp (A_t X_t \varepsilon_t + \frac{1}{2} A_t^2 X_t^2))} =
\frac{\hat{\pi}_t (J_t \exp (A_t X_t Y_t - \frac{1}{2} A_t^2 X_t^2))}{\hat{\pi}_t (\exp (A_t X_t Y_t - \frac{1}{2} A_t^2 X_t^2))},
$$
where $\hat{\pi}_t(\cdot)$ denotes a conditional expectation with respect to $\YY_{t}$ under $\hat{\Pg}$.
Note that under $\hat{\Pg}$ the distribution of $(X_s, Y_r)_{s\le t, \, r\le t-1}$  is the same as under $\Pg$ and
$Y_t$ is a ${\cal N}(0,1)$ random variable  independent of $(X_s, Y_r)_{s\le t, \, r\le t-1}$ .

To understand this point it is sufficient to write the following equality for the mutual characteristic function with arbitrary real numbers $(\alpha_{j}, \lambda_{j})$:
\begin{multline*}
\hat{\E} \exp \left\{i\sum_{j=1}^t \alpha_j X_j + i\sum_{j=1}^t \lambda_j Y_j \right\} =
\\
=\E \exp \left\{i\sum_{j=1}^{t} \alpha_j X_j + i\sum_{j=1}^{t-1} \lambda_j Y_j + i\lambda_{t} Y_t - A_t X_t \varepsilon_t - \frac{1}{2} A_t^2 X_t^2 \right\} = \\
= \E \left( \E \left. \exp \left\{i\sum_{j=1}^{t} \alpha_j X_j + i\sum_{j=1}^{t-1} \lambda_j Y_j + i\lambda_{t} Y_t - A_t X_t \varepsilon_t - \frac{1}{2} A_t^2 X_t^2 \right\} \right/ {\cal X }_{t}\right) =
\\
=\E  \exp \left\{i\sum_{j=1}^{t} \alpha_j X_j + i\sum_{j=1}^{t-1} \lambda_j Y_j + i\lambda_{t} A_t X_t - \frac{1}{2} A_t^2 X_t^2 + \frac{1}{2} (i\lambda_{t} - A_t X_t)^2 \right\} =
\\
= e^{-\frac{1}{2}\lambda_{t}^2} \E\exp \left\{i\sum_{j=1}^{t} \alpha_j X_j + i\sum_{j=1}^{t-1} \lambda_j Y_j \right\},
\end{multline*}
where ${\cal
X}_t$ is the  $\sigma$-field
${\cal
X}_t=\sigma(\{X_s\, ,\, 1\leq s\leq t\})$.
Hence,
\begin{multline*}
\hat{\pi}_t (J_t \exp(A_t X_t Y_t - \frac{1}{2}A_t^2 X_t^2 )) =
\\
=\pi_{t-1} (J_t \exp (A_t X_t y - \frac{1}{2} A_t^2 X_t^2 ))|_{y=Y_t} =
\\
= \pi_{t-1} (J_t \beta_t^y)|_{y=Y_t}.
\end{multline*}

Similarly,
$$
\hat{\pi}_t \left(\exp (A_t X_t y - \frac{1}{2} A_t^2 X_t^2 )\right) = \left.\pi_{t-1} (\beta_t^y)\right|_{y=Y_t}\,,
$$
and hence \eqref{pitJ} holds.

\textbf{II.} In the second step we will calculate the ratio $\frac{\mI_t}{\mI_{t-1}}$ which, due to \eqref{pitJ} can be rewritten as
\begin{equation}\label{ratio}
\frac{\mI_t}{\mI_{t-1}} = \frac{\pi_t(J_t)}{\pi_{t-1}(J_{t-1})}=\left. \frac{\pi_{t-1}(J_t \beta_t^y)}{\pi_{t-1}(J_{t-1})\pi_{t-1}(\beta_t^y)} \right|_{y=Y_t}.
\end{equation}

For this aim similarly to what we proposed in~\cite{mkalbmcr3} and \cite{AMM} we  introduce the auxiliary processes $(Y_t^2)_{t\ge 1}$ and $(\xi_t)_{t\ge 1}$.
Let $\bar{\varepsilon} =(\bar{\varepsilon}_t)_{t\ge 1}$ be a  sequence of i.i.d. $\mN(0,1)$ random variables independent of $X$ and define $(Y_t^2, \xi_t)_{t\ge 1}$ by:
\begin{equation}\label{Yaux}
Y_t^2= Q_t(X_t-h_t) + \sqrt{Q_t}\bar{\varepsilon}_t,
\end{equation}
\begin{equation}\label{xi:eq}
\xi_t=\sum_{s=1}^{t}(X_s-h_s) Y_s^2.
\end{equation}

Now the following equality holds:
$$
\left.\frac{\pi_{t-1}(J_t \beta_t^y)}{\pi_{t-1}(J_{t-1})}\right|_{y=Y_t} = \left.\frac{\bp_{t-1}(\exp \{-\frac{1}{2}Q_t (X_t-h_t)^2 - \xi_{t-1}\} \beta_t^y)}{\bp_{t-1}(\exp(-\xi_{t-1}))}\right|_{y=Y_t},
$$
where $\bp_t(\cdot)$ stands for a conditional expectation w.r.t. to the $\sigma$-field $\bar{\cal Y}_t=\sigma(\{Y_s, Y_s^2, {s\le t}\})$ under the initial measure $\Pg$.

Again the proof of this equality  is based on the Bayes formula. Namely, let $\tilde{\Pg}$ be the new probability measure defined by
\begin{equation}\label{ptilde}
\frac{d\tilde{\Pg}}{d\Pg} = \rho_{t-1}=\exp\left\{ -\frac{1}{2} \sum_1^{t-1} Q_s (X_s-h_s)^2 -\sum_1^{t-1} \sqrt{Q_s} (X_s-h_s) \bar{\varepsilon}_s \right\}.
\end{equation}

Then $J_t \rho_{t-1}= \exp \{-\xi_{t-1} -\frac{1}{2} Q_t (X_t-h_t)^2 \}$ and $J_{t-1} \rho_{t-1}= \exp \{-\xi_{t-1}\}$. Thus
\begin{multline*}
\left.\frac{\bp_{t-1} (\exp(-\xi_t - \frac{1}{2}Q_t (X_t-h_t)^2 ) \beta_t^y)}{\bp_{t-1} (\exp(-\xi_{t-1}))}\right|_{y=Y_t}=
\\
= \left.\frac{\E (J_t \beta_t^y \rho_{t-1}/\bar{\YY}_{t-1})}{\E (\rho_{t-1}/\bar{\YY}_{t-1})}\cdot
\frac{\E (\rho_{t-1}/\bar{\YY}_{t-1})}{\E \exp(J_{t-1}\rho_{t-1})/\bar{\YY}_{t-1})}\right|_{y=Y_t}=
\\
= \left.\frac{\tilde{\E}(J_t \beta_t^y/\bar{\YY}_{t-1})}{\tilde{\E}(J_{t-1}/\bar{\YY}_{t-1})}\right|_{y=Y_t} =\left. \frac{\pi_{t-1} (J_t \beta_t^y)}{\pi_{t-1}(J_{t-1})}\right|_{y=Y_t},
\end{multline*}
where the last equality holds because  under the probability measure $\tilde\Pg$ the distribution of  $(X_s, Y_s)_{s\le t}$
is the same as under the initial measure $\Pg$ and
 $(X_s, Y_s)_{s\le t-1}$  is independent of $(Y_{s}^{2})_{s\leq t-1}$.


Finally we have proved the following:
\begin{equation}\label{derivlogar}
\frac{\pi_t(J_t)}{\pi_{t-1}(J_{t-1})} = \left.\frac{\bp_{t-1}(\exp\left[-\xi_{t-1}+A_t X_t y -\frac{1}{2} Q_t (X_t -h_t)^2 - \frac{1}{2}A_t^2 X_t^2\right])}{\bp_{t-1}(\exp(-\xi_{t-1})) \pi_{t-1}(\beta_{t}^{y})}\right|_{y=Y_t}.
\end{equation}
 At this point we will use the conditionally Gaussian properties of $(X_{t},\xi_{t-1})$ w.r.t. $\bar{\YY}_{t-1}$ and Lemma 11.6 \cite{lipshi1} which says that for a Gaussian pair $(U,V)$  with mean values $m_{_{U}},m_{_{V}}$,  variances $\gamma_{_{U}},\gamma_{_{V}}$ and covariance $\gamma_{_{UV}}$

\begin{multline}\label{feqg}
\E\exp\left\{-\frac{1}{2} DU^2 + \lambda_1 U - \lambda_2 V \right\} = (1+D\gamma_{_{U}})^{-1/2} \times
\\
\times \exp\left\{ -\lambda_2 m_{V} + \frac{\lambda_2^2}{2} \gamma_{_{V}} - \frac{1}{2} \cdot\frac{D}{1+D\gamma_{_{U}}} (m_{_{U}}-\lambda_2 \gamma_{_{UV}})^2 \right. +
\\
+\left. \frac{\lambda_1^2 \gamma_{_{U}} + 2\lambda_1 (m_{_{U}}-\lambda_2\gamma_{_{UV}})}{2(1+D\gamma_{_{U}})} \right\},
\end{multline}
for any real numbers $\lambda_{1},\lambda_{2}$ and $D\ge 0$.
Indeed, in \eqref{derivlogar} we will apply this formula to $(U,V)=(X_{t},\xi _{t-1})$ given $\bar{\YY} _{t-1}$  with
$$
D=S_{t}=Q_t+A_t^2, \quad \lambda_2=1, \quad \lambda_1=A_t y+ Q_t h_t,
$$
in the numerator and
$D=\lambda_1=0 ,\quad \lambda_2=1 $ in the first factor of the denominator and again  to $(U,V)=(X_{t},\xi _{t-1})$ given ${\YY} _{t-1}$ with
$$
D=A_t^2, \quad \lambda_2=0, \quad \lambda_1=A_t y,
$$
in the second factor of the denominator.

Collecting the terms as coefficients for $h^{2}_{t}$ and $h_{t}$,  we obtain that
\begin{multline*}
\frac{\mI_t}{\mI_{t-1}} = \frac{(1+S_{t}\bg_t)^{-1/2}}{(1+A_{t}^{2}\gamma_t)^{-1/2}} \cdot \exp \left\{-\frac{Q_{t}}{2} \frac{1+A_t^2 \bg_t}{1+ S_t \bg_t} \times \left[ h_{t}-\frac{Z_{t}^h+A_t \bg_t Y_t}{1+A_t^2 \bg_t} \right]^2 \right\} \times \\
\times \exp\left\{
-\frac{A_t^2 (Z_{t}^h)^2 - A_{t}^{2} \bg_{t} Y_t^2}{2(1+A_t^2 \bg_{t})} + \frac{Y_tZ_{t}^hA_t}{1+A_{t}^2 \bg_t} + \frac{1}{2} \cdot \frac{A_t^2 \pi_{t-1}^2(X_{t}) - 2 A_t \pi_{t-1}(X_{t}) Y_t - A_{t}^2 Y_t^2 \gamma_t}{1+A_{t}^2\gamma_t}
\right\},
\end{multline*}
where
$Z_{t}^h=\bp_{t-1}(X_{t})-\bg_{_{X\xi}}(t)$ with
\begin{multline}\label{def gamma Xxi}
\bg_{_{X\xi}}(t)=\Eg[(X_t-\bp_{t-1}(X_t))(\xi_{t-1}-\bp_{t-1}(\xi_{t-1}))/
{\bar{\YY}}_{t-1}],\,t\ge 2\,;\,
\bg_{_{X\xi}}(1)=0\,.
\end{multline}

To finish the proof we just replace $Y_t$ by $\nu_t+A_t \pi_{t-1}(X_{t})$. Thus in the last exponential  term we find:
\begin{multline*}
\exp\left\{ - \frac{\nu_t^2 A_t^2 (\gamma_t-\bg_t)} {2(1+A_t^2 \bg_t)(1+A_t^2 \gamma_t)} + \frac{Z_{t}^{h}-\pi_{t-1}(X_{t})}{1+A_{t}^2\bg_t}A_t \nu_t \right. -
\\ - \left.\frac{1}{2}\cdot \frac{A_t^2}{1+A_t^2 \bg_t} (Z_{t}^h-\pi_{t-1}(X_{t}))^2 \right\},
\end{multline*}
which gives the Proposition.

\begin{rem}\label{probinterp}
\begin{enumerate}
  \item
    Note that now the probabilistic interpretation  of the ingredients $\bg_t$ and $Z_t^{h}$ is clarified for \textbf{negative} $\mu$. Namely,
$\bg_t=\E (X_t-\bp_{t-1}(X_t))^2,$
and $Z_t^{h}=\bp_{t-1}(X_t)-\bg_{_{X\xi}}(t)$, but when
   $\mu$ is \textbf{positive},  there is no such connection  anymore.
\item Observe that actually $\bp_{t-1}(X_t)$ and $\bg_{_{X\xi}}(t)$ are $\bar{\YY}_{t-1}$-measurable, but the difference $Z_t^{h}=\bp_{t-1}(X_t)-\bg_{_{X\xi}}(t)$ is $\YY_{t-1}$ measurable.
\end{enumerate}
\end{rem}

\paragraph{Proof of Theorem \ref{LEGsol}} The statement of Theorem
\ref{LEGsol} is the direct consequence of
Proposition \ref{CMC}. Indeed, we claim that the following chain of inequalities holds for any $h:\,h_t\in \YY_{t}, t\ge 1$ :
$$
\E \mu \exp \left\{\frac{\mu}{2} \sum_{t=1}^T (X_t-h_t)^2 Q_t \right\}
$$
$$
=\E\left[\E\mu \left(\left.\exp\left\{\frac{\mu}{2} \sum_{t=1}^T (X_t-h_t)^2 Q_t \right\} \right/ {\YY}_{T}\right)\right]
$$
$$
=\mu\E\prod_{t=1}^T \left[\frac{1+S_t \bg_t}{1+A_t^{2} \bg_t}\right]^{-1/2}  \times \exp\left\{\frac{\mu}{2} Q_t \frac{1+A_t^2 \bg_t}{1+ S_t \bg_t} \times \left[ h_t - \frac{Z_{t}^h+A_t \bg_t Y_t}{1+ A_t^2 \bg_t } \right]^{2}\right\} \times \mathcal{M}_T,
$$

$$
\stackrel{(a)}{\ge}\prod_{t=1}^T \left[\frac{1+S_t \bg_t}{1+A_t^{2} \bg_t}\right]^{-1/2} \mu \E \mathcal{M}_T
 $$

$$
\stackrel{(b)}{=} \mu \prod_{t=1}^T \left[\frac{1+S_t \bg_t}{1+A_t^{2} \bg_t}\right]^{-1/2} .
$$
Of course under condition $(C_{\mu})$, since the term in the last line is finite, it is sufficient to consider the case:
\begin{equation}\label{expfin}
\E \mu \exp \left\{\frac{\mu}{2} \sum_{t=1}^T (X_t-h_t)^2 Q_t \right\}
 < \infty,
\end{equation}
which gives the first equality. Inequality $(a)$ follows directly from Proposition~\ref{CMC}. Equality $(b)$ is a direct consequence of \eqref{feqg} which gives that $\E \mathcal{M}_T=1.$
Now, to obtain the lower bound
we must take
$$
\bh_{t}= \displaystyle{\frac{Z_{t}^{\bar{h}}+A_t \bg_t Y_t}{1+ A_t^2 \bg_t}},\, t\ge 1,
$$
or equivalently

$$
\bh_{t}= Z_{t}^{\bh} +  \frac{A_t \bg_{t}}{1+A_t^2 \bg_{t}} (Y_t - A_t Z_{t}^{\bh}),\, t\ge 1,
$$
where $Z^{h}$ is the solution
of equation (\ref{eqzh}), which means that
$$
Z_{t}^{\bh} =m_t +\sum_{l=1}^{t-1} \frac{\bg(t,l) A_l}{1+A_l^2 \bg_l}  [Y_l -A_l Z_{l}^{\bh}],
$$
and hence
$$
\bh_{t}= m_t + \sum_{l=1}^{t} \frac{\bg(t,l) A_l}{1+A_l^2 \bg_l}  [Y_l -A_l Z_{l}^{\bh}]=m_t+\sum_{l=1}^{t} A_l \bg(t,l)  (Y_l -A_l \bh_l ).
$$
Thus
$\bar{h}$ is the unique solution of equation \eqref{hbar}.
Finally for $\bar{h}$ the lower bound is attained.
\begin{rem}\label{probinterp'}
\begin{enumerate}
   \item   It is worth  emphasizing that the process $\displaystyle{\tilde{Z}^{h}_{t}=\frac{Z_{t}^h+A_t \bg_t Y_t}{1+ A_t^2 \bg_t }}$  is the solution of the following recursive equation:

 \begin{equation}\label{zhtild}
 \tilde{Z}^{h}_{t}=m_t -\sum_{l=1}^{t-1}\bg(t,l) \frac{\mu Q_l}{1+S_l \bg_l} (h_l - \tilde{Z}_{l}^h) + \sum_{l=1}^{t}\bg(t,l) A_l (Y_l - A_l \tilde{Z}_{l}^{h}),
 \end{equation}
 and hence the equality $\bh_t =\displaystyle{\widetilde{Z}^{\bh}_{t}} $ implies immediately the equation \eqref{eqzh} for $\bh$.
 This process  $\tilde{Z}^{h}$ also has
  a probabilistic interpretation as well as $\displaystyle{\tilde{\gamma}_{t}=\frac{ \bg_t}{1+ A_t^2 \bg_t }}$. This interpretation will be given in Section~\ref{interpret}.
\end{enumerate}
\end{rem}

\subsection{Solution of the auxiliary filtering  problems}\label{AFP}
Here, for an arbitrary Gaussian sequence $X$, we deal with the
one-step prediction
and filtering problems of the signals $X$
and $\xi$ given by \eqref{xi:eq} respectively from the observation of
$\bar{Y}=(Y,Y^{2})$ defined in \eqref{observ} and \eqref{Yaux}.
Actually, we follow the ideas proposed in our paper
\cite{mkalbmcr3}.
Recall that the solutions can
be reduced to equations for the conditional moments. The following
statement provides the equations for the characteristics which give
the solution
of the prediction problem and the equation for the other quantity
$\bp_{t-1}(X_t)-\bg_{_{X\xi}}(t)$ appearing  in Proposition~\ref{p:CM} for $\mu=-1.$
\begin{theorem}\label{filter}
The conditional mean $\bp_{t-1}(X_t)$ and the
variance of the one-step prediction error
$\bg_t=\Eg[X_t-\bp_{t-1}(X_t)]^{2}$ are given by the equations
\begin{multline}\label{eq1 filter}
\bp_{t-1}(X_t)=m_{t}+\sum\limits_{s=1}^{t-1}
\frac{\bg(t,s)}{1+(A_s^2+Q_s)\bg_s}[A_s(Y_{s}-A_s\bp_{s-1}(X_s)) \\
+Q_s(Y_s^2 - Q_s(\bp_{s-1}(X_s)-h_s)]\,,\quad t\geq 1,
\end{multline}
\begin{equation}\label{eq2 filter}
\bg_t=\bg(t,t)\,,\quad t\geq 1\,.
\end{equation}
where $\bg=(\bg(t,s),\, 1\le s \le t)$ is the unique solution of equation \eqref{GAMMABAR}.
Moreover, with $\bg_{_{X\xi}}(t)$ defined by (\ref{def gamma Xxi}), the difference $\bp_{t-1}(X_t)-\bg_{_{X\xi}}(t)$ is the
solution $Z^h_t$ of equation (\ref{eqzh}).
\end{theorem}
\paragraph{Proof} Note that since $h_t \in \YY_t$  and the joint
distribution of $(X_{r},Y_{s}, Y_s^2+Q_s h_s)$ for any $r\,,s$  is Gaussian, we
can apply the Note following
  Theorem 13.1 in \cite{lipshi1bis}. For any
$k\le t$ we can write
\begin{equation}\label{piX}
\left\{
\begin{array}{l}
\bp_k(X_t)=\bp_{k-1}(X_t) + [\cov(X_t,\overline{\nu}_k)]^{\prime} \var(\overline{\nu}_k)^{-1} \overline{\nu}_k,
\\
\bp_{0}(X_{t})=m_{t},
\end{array}
\right.
\end{equation}
where
$$
\overline{\nu}_{k}=\bar{Y}_{k}-\Eg(\bar{Y}_{k}/\overline{\YY}_{k-1})=\left(
                               \begin{array}{c}
                                Y_k-A_k\bp_{k-1}(X_k) \\
                                Y_k^2 +Q_kh_k - Q_k\bp_{k-1}(X_k)
                               \end{array}
                             \right)
 $$
is the innovation with covariance matrices
\begin{equation}\label{varnu}
\var(\overline{\nu}_k)=\left(
                               \begin{array}{cc}
                                1+A_k^2\bg_k & A_kQ_k\bg_k\\
                                A_kQ_k\bg_k & Q_k+Q_k^2\bg_k
                               \end{array}
                             \right),
\end{equation}
and
\begin{equation}\label{covarXnu}
\cov(X_t,\overline{\nu}_k) = \bg (t,k) \left( \begin{array}{c} \! A_k \! \\ \! Q_k \! \end{array} \right),
\end{equation}
 with
\begin{equation}\label{def gamma t s}
\bg(t,k)=\Eg (X_t-\bp_{k-1}(X_t))(X_{k}-\bp_{k-1}(X_k))\,.
\end{equation}
By the definition (\ref{def gamma t s}), we see for $k=t$ that the
variance $\bg_t$ is given by (\ref{eq2 filter}). Now,
equality (\ref{piX}) implies
\begin{multline}\label{pikt}
\bp_k(X_t)=m_t+\sum_{l=1}^k \bg(t,l) \left( \begin{array}{cc} \! A_l \! & \! Q_l \! \end{array} \right)
(\var \overline{\nu}_l)^{-1} \overline{\nu}_l=
\\
=m_{t}+\sum\limits_{s=1}^{k}
\frac{\bg(t,s)}{1+(A_s^2+Q_s)\bg_s}[A_s(Y_{s}-A_s\bp_{s-1}(X_s))+
\\
+Q_s(Y_s^2 - Q_s(\bp_{s-1}(X_s)-h_s)],
\end{multline}
and putting $k=t-1$ we get nothing but equation (\ref{eq1
filter}). Concerning the solution of the one-step prediction
problem, it just remains to show that the covariance $\bg(t,s)$
satisfies  equation (\ref{GAMMABAR}).

Let us define
$$
\delta_{X} (t,l) = X_t - \bp_l (X_t)\,.
$$
According to (\ref{piX}) we can write
$$
\delta_X (t,l) = \delta_X (t,l-1) - \bg(t,l)\left( \begin{array}{cc} \! A_l \! & \! Q_l \! \end{array} \right)
(\var \overline{\nu}_l)^{-1} \overline{\nu}_l,
$$
and so

\begin{multline*}
\E \delta_X(t_1, l) \delta_X(t_2,l) = \E \delta_X(t_1, l-1) \delta_X(t_2, l-1) -
\\
-\bg(t_1,l) \bg(t_2,l) \left({A_l \atop Q_l}\right)^{\prime} \var(\bar{\nu}_l)^{-1} \left({A_l \atop Q_l}\right),
\end{multline*}
or
\begin{multline}\label{deltadelta}
\Eg \delta_{X}(t^1,l)\delta_{X}(t^2,l)=\Eg
\delta_{X}(t^1,0)\delta_{X}(t^2,0)- \\
-\sum\limits_{r=1}^{l}\bg(t,r) \bg(s,r) \frac{A_r^2+Q_r}{1+ (A_r^2 +Q_r ) \bg_{r}}.
\end{multline}
Taking $t^1=t\,, t^2=s\,, l= s-1$ in (\ref{deltadelta}), it is
readily seen that
equation (\ref{GAMMABAR}) holds for $\bg(t,s)$.

Now we analyze the difference $\bp_{t-1}(X_t)-\bg_{_{X\xi}}(t)$.
Using the representation $\xi_t=\sum_{s=1}^t (X_s-h_s) Y_s^2$
we can rewrite $\bp_{t-1}(\xi_{t-1})$ in the following form

$$
\bp_{t-1}(\xi_{t-1})= \sum_{s=1}^{t-1} (\pi_{t-1}(X_s)-h_s) Y_s^2,
$$
which implies that
$$
\xi_{t-1}-\bp_{t-1}(\xi_{t-1})= \sum_{s=1}^{t-1} (X_s - \bp_{t-1}(X_s)) Y_s^2.
$$
So we have
\begin{multline}\label{ga X xi}
\bg_{_{X\xi}}(t) = \sum_{s=1}^{t-1} \bp_{t-1} [(X_s-\bp_{t-1}(X_s)) (X_t-\bp_{t-1}(X_t))] Y_s^2 =
\\
= \sum_{s=1}^{t-1} \E (X_s-\bp_{t-1}(X_s)) (X_t-\bp_{t-1}(X_t)) Y_s^2
= \sum_{s=1}^{t-1} \tg(t,s) Y_s^2,
\end{multline}
where
\begin{equation}\label{ad gamma }
\tg(t,s)=\E (X_s-\bp_{t-1}(X_s)) (X_t-\bp_{t-1}(X_t))= \bg(s,t).
\end{equation}
Using the definitions (\ref{def gamma t s}) and (\ref{ad gamma })
we can write
$$
\tg(t,s)-\bg(t,s) = - \E X_t (\bp_{t-1}(X_s) - \bp_{s-1}(X_s)).
$$

Again, applying  the Note following Theorem 13.1 in \cite{lipshi1bis}, we
can write also
$$
\bp_{l}(X_{r})= \bp_{l-1}(X_{r})+\bg(t,l) \left( \begin{array}{cc} \! A_l \! & \! Q_l \! \end{array} \right)
(\var \overline{\nu}_l)^{-1} \overline{\nu}_l.
$$
This means that
$$
\pi_{t-1}(X_{r})- \pi_{r-1}(X_{r})= \sum\limits_{l=r}^{t-1}\bg(t,l) \left( \begin{array}{cc} \! A_l \! & \! Q_l \! \end{array} \right)
(\var \overline{\nu}_l)^{-1} \overline{\nu}_l\,,
$$
or equivalently
$$
\pi_{t-1}(X_{r})- \pi_{r-1}(X_{r})
=\sum\limits_{l=r}^{t-1}\tg(l,t) \left( \begin{array}{cc} \! A_l \! & \! Q_l \! \end{array} \right) (\var \overline{\nu}_l)^{-1}
\overline{\nu}_l\,.
$$
Then, multiplying by $X_t$ and taking expectations in both sides, we get
\begin{multline*}
\Eg X_{t}
(\pi_{t-1}(X_{r})-\pi_{r-1}(X_{r}))
=\sum\limits_{l=r}^{t-1}\tg(l,r) \left( \begin{array}{cc} \! A_l \! & \! Q_l \! \end{array} \right) (\var \overline{\nu}_l)^{-1}\cov(X_t,\overline{\nu}_l)=
\\
=\sum_{l=s}^{t-1} \tg(l,s) \bg(t,l)\frac{A_l^2 +Q_l}{1+(A_l^2 +Q_l)\bg_l} .
\end{multline*}
Hence we have proved the following relation
\begin{equation}\label{dif gamma}
\tg(t,s)-\bg(t,s) =-\sum_{l=s}^{t-1} \tg(l,s) \bg(t,l)\frac{A_l^2 +Q_l}{1+(A_l^2 +Q_l)\bg_l} .
\end{equation}

Now we can show that the difference $Z^h_{t}=\bp_{t-1}(X_t)-\bg_{_{X\xi}}(t)$
satisfies the equation (\ref{eqzh}). Using (\ref{pikt}) and (\ref{ga X xi}),
we can write

\begin{multline}\label{zhbeg}
Z^h_{t} = m_t+\sum_{l=1}^{t-1} \bg(t,l) \left( \begin{array}{cc} \! A_l \! & \! Q_l \! \end{array} \right) (\var \overline{\nu}_l)^{-1} \overline{\nu}_l - \sum_{s=1}^{t-1} \tg(t,s) Y_s^2 =
\\
= m_t +\sum_{l=1}^{t-1} \frac{A_l\bg(t,l)} {1+(A_l^2 +Q_l)\bg_l}(Y_l - A_l \bp_{l-1}(X_l)) +
\\
+
\sum_{l=1}^{t-1} \frac{\bg(t,l)} {1+(A_l^2 +Q_l)\bg_l}(Y_l^2-Q_l (\bp_{l-1}(X_l)-h_l))
-\sum_{l=1}^{t-1} \tg(t,l) Y_l^2 =
\\
= m_t +\sum_{l=1}^{t-1}\frac{A_l\bg(t,l)} {1+(A_l^2 +Q_l)\bg_l} Y_l + \sum_{l=1}^{t-1} \frac{\bg(t,l)} {1+(A_l^2 +Q_l)\bg_l} Q_l h_l - \\
 - \sum_{l=1}^{t-1} \bg(t,l) \frac{A_l^2 +Q_l}{1+(A_l^2 +Q_l)\bg_l} \bp_{l-1}(X_l) +
\\
+ \sum_{l=1}^{t-1} [\frac{\bg(t,l)} {1+(A_l^2 +Q_l)\bg_l} - \tg (t,l)] Y_l^2.
\end{multline}
Now we can rewrite the last term in \eqref{zhbeg} using the equality \eqref{dif gamma}. We have
\begin{multline}\label{gammause}
\sum_{l=1}^{t-1} [\frac{\bg(t,l)} {1+(A_l^2 +Q_l)\bg_l} - \tg (t,l)] Y_l^2= \sum_{l=1}^{t-1} \bg(t,l) (\frac{1} {1+(A_l^2 +Q_l)\bg_l}-1) Y_l^2 +
\\
+ \sum_{l=1}^{t-1} \sum_{r=l}^{t-1} \bg(t,r) \tg(r,l)\frac{A_r^2 +Q_r}{1+(A_r^2 +Q_r)\bg_r} Y_l^2 = \\
= \sum_{r=1}^{t-1} \bg(t,r)
 \left[\sum_{l=1}^{r-1} \tg(r,l) Y_l^2 \right]
\frac{A_r^2 +Q_r}{1+(A_r^2 +Q_r)\bg_r}=
\\
=\sum_{r=1}^{t-1} \bg(t,r)\bg_{_{X\xi}}(r)\frac{A_r^2 +Q_r}{1+(A_r^2 +Q_r)\bg_r}\,,
\end{multline}
where in the last step we have used equality (\ref{ga X xi}).

Finally \eqref{zhbeg}-\eqref{gammause} imply:
\begin{multline*}
Z^h_{t} = m_t +\sum_{l=1}^{t-1}\frac{A_l\bg(t,l)} {1+(A_l^2 +Q_l)\bg_l} Y_l + \sum_{l=1}^{t-1} \frac{Q_l\bg(t,l)} {1+(A_l^2 +Q_l)\bg_l} h_l -
\\
 - \sum_{l=1}^{t-1} \bg(t,l) \frac{A_l^2 +Q_l}{1+(A_l^2 +Q_l)\bg_l}[ \bp_{l-1}(X_l)-\bg_{_{X\xi}}(l)]=
\\
 m_t +\sum_{l=1}^{t-1}\frac{A_l\bg(t,l)} {1+(A_l^2 +Q_l)\bg_l} Y_l + \sum_{l=1}^{t-1} \frac{Q_l\bg(t,l)} {1+(A_l^2 +Q_l)\bg_l} h_l -
\\
 - \sum_{l=1}^{t-1} \bg(t,l) \frac{A_l^2 +Q_l}{1+(A_l^2 +Q_l)\bg_l}Z^h_{l}\,,
\end{multline*}
which is nothing else but equation (\ref{eqzh}) with $\mu=-1$.

\section{Particular cases and applications}\label{PC}

Here we deal with some specific cases where the results can be further elaborated. For two examples we can apply directly Theorem~\ref{LEGsol} and moreover
the special structure of the covariances allows  to simplify the answer.
\subsection{LEG filtering of Gauss-Markov sequences }\label{GMS}
In this part we concentrate on the case of a
Gaussian AR(1) process $X$, {\em i.e.}, a Gauss-Markov process
driven by
\begin{equation}\label{model AR}
X_{t}= a_{t} X_{t-1}+D_t^\frac{1}{2}\widetilde{\varepsilon}_{t}\,,\; t\ge 1\,;
\quad X_{0}=x\,,
\end{equation}
where $(\widetilde\varepsilon_{t}, \; t=1,2,\dots)$ is a sequence of i.i.d.
standard Gaussian random variables and
$(D_t,\,t\ge 1)$ is a (deterministic)
sequence of real numbers such that $D_t\ge 0$ for $t\ge 1$. In this setting,
it is easy to check that the mean
and covariance functions of $X$ are given by
$$
m_{t}=[\prod_{u=1}^{t}a_u]x=\Lambda_{t}x\,;\quad
K(t,s)=[\prod_{u=s+1}^{t}a_u]k_s=\frac{\Lambda_{t}}{\Lambda_{s}}k_s\,,\;1\le s\le t\,,
$$
where $\Lambda_{t}=\prod\limits_{u=1}^{t}a_u$ and
$$
k_t=a_t^2k_{t-1} +D_t,\, t\ge 1,\,k_0=0.
$$
Suppose that the following the Riccati type  equation
\begin{equation}\label{ricmar}
\bg_s=D_s+\frac{a_{s}^2 \bg_{s-1}}{1+(A_{s-1}^2-\mu Q_{s-1})\bg_{s-1}},\, s\ge 1,\,\bg_{0}=0,
\end{equation}
has a unique nonnegative solution.

From the classical filtering  theory it is
well-known that (for $\mu<0$ ) $\bg_s$ is nothing but the variance  of the error  of the one-step prediction problem
 of the signal $X$ given by the auxiliary observation $\bar{Y}$ defined by equations \eqref{observ} and \eqref{Yaux}. Then, it is readily
seen that the function $\bg(t,s)$, where
$\displaystyle{\bg(t,s)=\frac{\Lambda_{t}}{\Lambda_{s}}\bg_s}$ is the
solution of equation (\ref{GAMMABAR}) and that moreover equation
\eqref{hbar} for the solution $\bh$ of the LEG filtering problem (\ref{LEGdef}) can be reduced to the
following one:
\begin{equation}\label{rsmkk}
\bh_t=\frac{a_t}{1+A_t^2\bg_t}\bh_{t-1} +\frac{A_t \bg_t}{1+A_t^2\bg_t}Y_t, \, t\ge 1,\, \bh_0=x,
\end{equation}
or, equivalently:
\begin{equation*}
\bh_t=a_t \bh_{t-1} +\frac{A_t \bg_t}{1+A_t^2\bg_t} [Y_t-a_tA_t \bh_{t-1}], \, t\ge 1,\, \bh_0=x.
\end{equation*}
Actually  equation \eqref{rsmkk} can also be obtained directly from the general filtering theory (for $\mu=-1$ and replacing $Q$ by $-\mu Q$).
For arbitrary $(h_{t} \in \YY_{t}, t\ge 1)$  the Note following
  Theorem 13.1 in \cite{lipshi1bis} gives the equation for $Z^{h}$:
\begin{multline*}
Z^{h}_{t}=a_t Z^{h}_{t-1} +a_t \bg_t\frac{Q_{t-1} }{1+S_{t-1}\bg_t} [h_{t-1}-Z^h_{t-1}] +
\\
+a_t \bg_t\frac{A_{t-1} }{1+S_{t-1}\bg_t} [Y_{t-1}-A_{t-1}Z^h_{t-1}], \, t\ge 1,\, Z^h_0=x.
\end{multline*}

Hence, again  the solution
$\displaystyle{\bh_{t}= \displaystyle{\frac{Z_{t}^{\bar{h}}+A_t \bg_t Y_t}{1+ A_t^2 \bg_t}},\, t\ge 1,
}$ of the LEG filtering problem (\ref{LEGdef}) is given by \eqref{rsmkk}.

Let us emphasize that these equations are nothing but those given in  Speyer  \textit{et al.}
\cite{speyer}.

It is interesting to note that in the case $a_t=0$ (i.i.d. signal) the solution of the LEG filtering problem is
nothing else but the solution of the risk neutral filtering problem \textit{i.e.} $\bh_t=\pi_t(X_t)$.


\subsection{LEG filtering of  moving averages of order 1}\label{MA1X}
Here we consider the case of a MA(1) process, {\em i.e}, a non
Markovian process $X$ defined by
$$
X_{t} = \widetilde{\varepsilon}_{t} + \lambda \widetilde{\varepsilon}_{t-1}\,; t\ge 1\,,
$$
where $( \widetilde{\varepsilon}_0, \widetilde{\varepsilon}_1,\dots)$~is a sequence of i.i.d. standard
Gaussian variables and $\lambda$ is a real number.  Of course
$X$ is centered and has the
covariance function
$K(t,s)= 1+\lambda^{2}$ if $s=t$, $\lambda$ if $s=t-1$ and
$0$ if $s < t-1$.
In order to solve equation (\ref{GAMMABAR}) we can take
$$
\bg(t,s) =0\,,\; s < t-1\,; \quad  \bg(t,t-1) = \lambda\,,\;
t\ge 1\,,
$$
and $\bg(t,t)=\bg_{t}$ where $\bg_{t}$  is
the solution of the equation:
$$
\bg_{t} =1+\lambda^{2}
 - \lambda \frac{A_{t-1}^2-\mu Q_{t-1}}{1+(A_{t-1}^2-\mu Q_{t-1})\bg_{t-1}}\,,\;t\ge 1\,;\quad \bg_0=1+\lambda^{2},
$$
provided that this equation has a unique nonnegative solution.

Moreover equation
\eqref{hbar} for the solution $\bh$ of the LEG filtering problem (\ref{LEGdef}) can be reduced to the
following one:
\begin{equation*}
\bh_t=\lambda\frac{A_{t-1}}{1+A_t^2\bg_t}[Y_{t-1}-A_{t-1}\bh_{t-1}] +\frac{A_t \bg_t}{1+A_t^2\bg_t}Y_t, \, t\ge 1,\, \bh_0=0.
\end{equation*}
Again, it is interesting to note that for $\lambda=0$ (i.i.d. signal) the solution of LEG filtering problem is
nothing else but the solution of the risk neutral filtering problem \textit{i.e.} $\bh_t=\pi_t(X_t)$.

\section{LEG and RS filtering problems}\label{disc}
Here, at first we show that actually the LEG and RS filtering problems have the same solution. Then we give an example which shows that in a more general context similar problems may have  different solutions.
\subsection{Equivalence of LEG and RS filtering problems}
Let $\bar{h}=(\bh_s)_{s \ge 1}$ be the solution of the LEG filtering problem (\ref{LEGdef})
given by equation \eqref{hbar}.
For any fixed $t \le T$, let us denote by $\hat{g}_{t}:$
\begin{multline*}
\hat{g}_{_t}=
\displaystyle{\arg\min_{g \in {\cal Y}_t} }\Eg\Big[\mu
\exp\Big\{\displaystyle{\frac{\mu}{2}(X_{t}-g)^{2}Q_{t}}
+\left.\displaystyle{\frac{\mu}{2} \sum_{s=1}^{t-1}
(X_{s}-\bar{h}(s))^{2}Q_{s} }\Big\} \right/
{\cal Y}_t\Big],
\end{multline*}
where  $g\in {\cal Y}_t$ means that $g$ is a ${\cal Y}_t$-measurable variable.
It follows directly  from
Proposition~\ref{p:CM}
 that, provided that $\displaystyle{1 + S_t\bg_t}> 0$, the equality
 $\displaystyle{\hat{g}_{t}= \displaystyle{\frac{Z_{t}^{\bar{h}}+A_t \bg_t Y_t}{1+ A_t^2 \bg_t}},\, t\ge 1}$ holds.
 Since it was noted in the proof of Theorem~\ref{LEGsol}  that $\displaystyle{\bh_{t}= \displaystyle{\frac{Z_{t}^{\bar{h}}+A_t \bg_t Y_t}{1+ A_t^2 \bg_t}},\, t\ge 1,}$
 hence we have also $\hat{g}_{t}=\bar{h}_t$. It means that for $ t\ge 1 $
the solution $\bar{h}$ of the LEG filtering problem  satisfies the following recursive
equation:
\begin{multline}\label{rde}
\hat{g}_{_t}=
\displaystyle{\arg\min_{g \in {\cal Y}_t} }\Eg\Big[\mu
\exp\Big\{\displaystyle{\frac{\mu}{2}(X_{t}-g)^{2}Q_{t}}
+\left.\displaystyle{\frac{\mu}{2} \sum_{s=1}^{t-1}
(X_{s}-\bar{h}(s))^{2}Q_{s} }\Big\} \right/
{\cal Y}_t\Big].
\end{multline}
Indeed, in the literature, the recursion  (\ref{rde}) is the basic
\textbf{definition}  of the so-called risk-sensitive (RS)
filtering problem which was introduced in \cite{elliott3}. Therefore we have also proved the following
statement

\begin{theorem}\label{RS}
Assume that the condition $(C_{\mu})$ is
satisfied.   Let $\bh=(\bh_t)_{t\ge 1}$  be the unique solution of  equation
(\ref{hbar}), \textit{i.e.}, $\bar{h}$  is the solution of the LEG
filtering problem \eqref{LEGdef}.
Then $\bar{h}$  is the solution of the RS filtering problem
\eqref{rde}.
\end{theorem}

\subsection{Discrepancy between LEG and RS type filtering problems: an example }\label{EX}
Actually, we did not find in the literature any trace of the
discussion
 about the relationship between the LEG
filtering problem (\ref{LEGdef}) and the RS filtering
problem (\ref{rde}) even in a Gauss-Markov  case. As a complement to our observation that these two problems
have the same solution, we propose  an example to show that in a bit more general setting,
two similar problems may have different solutions.

For  given  positive  symmetric deterministic $2\times 2$ matrices
$\Lambda_s,1\le s\le T$, let us set $\Phi_{t}(h)=(X_{t} \, h_t)\Lambda_{t}\left(%
\begin{array}{c}
  X_{t} \\
  h_t \\
\end{array}%
\right)$.  We can define $\bar{h}_{t}\in {\YY}_{t},\, t\ge 1$ as a solution of a
\textit{LEG type filtering problem} :

\begin{equation}\label{defrssex}
\bh= \arg\min_{h_{t}\in \YY_{t},\, t\ge 1}\Eg \left[\mu \exp\left\{\frac{\mu}{2}
\sum_1^T \Phi_{s}(h)\right\}\right].
\end{equation}

We can also define
$\hat{h}$  as the solution of
the following recursive equation (\textit{RS type filtering problem}):
\begin{equation}\label{rdeex}
\hat{h}_{t}= \displaystyle{\arg\min_{g \in {\cal Y}_t}
}\,\Eg\Big[\mu
\exp\Big\{\displaystyle{\frac{\mu}{2}\Phi_{t}(g)\,}+\left.\displaystyle{\frac{\mu}{2} \sum_{1}^{t-1}
\Phi_{s}(\hat{h})\, }\Big\} \right/ {\cal Y}_t\Big].
\end{equation}

The question which we discuss now is the following: does the equality
$\bar{h}=\hat{h}$ hold?

As we have just proved, the answer is positive for singular matrices
$\Lambda$, namely, when
$\Lambda_{11}=\Lambda_{22}=-\Lambda_{12}=Q$.
But in the general situation the answer may be negative. Actually it is sufficient to consider the following example:
 $\Lambda =\left(%
\begin{array}{cc}
  2 & -1 \\
  -1 & 1 \\
\end{array}%
\right),\, A_{t}=1,\,\mu=-1 $ and $X_{t}= X_{t-1}+\widetilde{\varepsilon}_{t}$, where $(\widetilde\varepsilon_{t}, \; t=1,2,\dots)$
is a sequence of i.i.d.
standard Gaussian random variables. Even in this Markov
case $\hat{h}\ne \bh $. More explicitly let us introduce the new probability measure $\hat{\Pg}:$
$$
\frac{d\hat{\Pg}}{d\Pg}= \frac{\exp\left[-\frac{1}{2}\sum\limits_{i=1}^{T}X_i^2\right]}{\Eg\exp\left[-\frac{1}{2}\sum\limits_{i=1}^{T}X_i^2\right]}.
$$
One can check that with respect to $\hat{\Pg}$  the observation model $(X_t,Y_t)_{t \ge 1}$ can be written in the following form:
$$
X_t= a_tX_{t-1} + D_t^\frac{1}{2}\hat{\varepsilon}_{t}\,,\; t\ge 1\,;
\quad X_{0}=x\,,
$$
$$
Y_t=X_t +\varepsilon_t\,,
$$
where $(\hat\varepsilon_{t})_{t\ge 1}$ is a sequence of i.i.d.
standard Gaussian random variables independent of the sequence $\varepsilon$,
$$
a_t=D_t=\frac{1}{1+ \Gamma(T,t)}\,,
$$
and  $\Gamma(T,\cdot)$ is the solution of the backward Riccati equation
$$
\Gamma(T,t) = 1 +\frac{\Gamma(T,t+1)}{1+\Gamma(T,t+1)},\,\Gamma(T,T)=0.
$$
It can be checked that
$$
\Gamma(T,t)=10\frac{\lambda^{T}-\lambda^{t}}{(1-\sqrt{5})\lambda^{T}-(1+\sqrt{5})\lambda^{t}},\, \lambda = \frac{(3-\sqrt{5})}{(3+\sqrt{5})}.
$$
Indeed to explain this change of the observation model it is sufficient to calculate the conditional characteristic function:
$$
\hat{\Eg}\left[\left.\exp (i\lambda X_t)\right/{\cal
X}_{t-1}\right]=\frac{\Eg\left[\left.\exp\left[i\lambda X_t-\frac{1}{2}\sum\limits_{i=1}^{T}X_i^2\right] \right/{\cal
X}_{t-1}\right]}{\Eg\left[\left.\exp\left[-\frac{1}{2}\sum\limits_{i=1}^{T}X_i^2\right] \right/{\cal
X}_{t-1}\right]},
$$
where ${\cal
X}_{t-1}$ is the
$\sigma$-field
${\cal
X}_{t-1}=\sigma(\{X_s\, ,\, 1\leq s\leq t-1\})$.
But it follows directly from the equation (19)-(20) in \cite{mkalbmcr3} and from \eqref{feqg}
that
$$\hat{\Eg}\left[\left.\exp (i\lambda X_t)\right/{\cal
X}_{t-1}\right]=\exp\left\{\frac{i\lambda}{1+ \Gamma(T,t)}X_{t-1} -\frac{\lambda^{2}}{2(1+ \Gamma(T,t))}\right\}.
$$

Since the density $\displaystyle{\frac{d\widehat{\Pg}}{d\Pg}}$ does not depend on $h$, the initial  LEG filtering problem \eqref{defrssex} can be rewritten as:
$$
\bh= \arg\min_{h_{t}\in \YY_{t},\, 1\le t\le T}\hat{\Eg} \left[- \exp\left\{-\frac{1}{2}
\sum_1^T (X_s-h_s)^{2}\right\}\right].
$$
Hence we can apply Theorem~\ref{LEGsol} or in particular \eqref{ricmar} and \eqref{rsmkk}. Clearly, $\bh$ depends on $T$ and $\hat{h}$ does not depend on $T$ by the definition. A bit more explicitly we have for example for $t=1$:
$
\displaystyle{\bh_{1}=\frac{1+\Gamma(T,1)}{2+\Gamma(T,1)}Y_1}
$
and obviously
$
\displaystyle{\hat{h}_{1}=\frac{\pi_{1}(X_1)}{1+\gamma_{1}}=\frac{1}{4}Y_1}
$
and clearly they are different.

\section{Information state, interpretation}\label{interpret}
In this section we discuss the probabilistic interpretation of the ingredients of the ``information state'' which was introduced
in the context of RS filtering and LEG control problems.
By the definition, the ``information state'' contains all the information needed to describe the solution of the concerned optimization problem. In particular it takes into account  the cost function but not only estimates of the signal and it should
give the total information about the model states available in the measurement.

\medskip
\textit{Risk-Sensitive Filtering}

\medskip\noindent
In the context of the RS filtering problem  the definition of the information state can be found for example in  \cite{elliott1}.
It is the density $\lambda_{t}$, with respect to the Lebesgue measure, of the non normalized random measure
$\omega_{t}$:
\begin{equation}\label{mescond}
 \omega_{t}(dx)=\displaystyle{\Eg\left[\left.\mathbb{I}(X_t \in dx)\exp\left\{\frac{\mu}{2} \sum_{s=1}^{t-1}
(X_{s}-h(s))^{2}Q_{s}\right\} \right/ {\cal
Y}_t\right]},
\end{equation}
 where $h_t\in \YY_t, t\ge 1 $
and  the observation $Y$ is defined by the equation \eqref{observ}.\\
In a classical Gauss-Markov setting, an explicit representation of $\lambda_{t}$ can be obtained as the  solution of some
recurrence  equation (see, \textit{e.g.},
\cite{collings}).\\
We claim that for a general Gaussian signal $X$ the density $\lambda_{t}$
 satisfies the following equality:
\begin{multline}\label{cdensg}
\lambda_{t}(x) =
\displaystyle{\frac{1}{\sqrt{2\pi \tg_t}}
\exp\left\{-\frac{(x-\widetilde{Z}_{t}^h)^{2}}{2\widetilde{\gamma}_{t}}
 \right\}}\times
\\
\prod_{r=1}^{t-1} \left[\frac{1+S_r \bg_r}{1+A_r^{2} \bg_r}\right]^{-1/2} \times \exp\left\{\frac{\mu}{2} Q_r \frac{1+A_r^2 \bg_r}{1+ S_r \bg_r} \times \left[ h_r -\widetilde{Z}_{t}^h \right]^{2}\right\} \times \mathcal{M}_t,
\end{multline}
where $\displaystyle{ \widetilde{Z}_{t}^h=\frac{Z_{t}^h+A_t \bg_t Y_t}{1+A_t^2\bg_t}}$ is the solution of the equation \eqref{zhtild}, $\displaystyle{ \widetilde{\gamma}_{t}= \frac{\bg_{t}}{1+A_t^2\bg_t}}$,  $\bg, Z^h$ are the solutions of
equations (\ref{GAMMABAR}) and (\ref{eqzh}) respectively and the martingale  $(\mathcal{M}_t)_{\ge 1}$ is defined by \eqref{martin}.

Indeed, to prove (\ref{cdensg})  it is sufficient to write the following:
\begin{equation}\label{bayesford}
\omega_{t}(dx) = \frac{\Eg[\mathbb{I}(X_t \in dx)\exp(-\xi_{t-1})/
{\cal \overline{Y}}_{t,t-1}]}{\Eg[\exp(-\xi_{t-1})/
{\cal \overline{Y}}_{t,t-1}]}\Eg\left[\left.\exp\left\{\frac{\mu}{2} \sum_{s=1}^{t-1}
(X_{s}-h(s))^{2}Q_{s}\right\}\right/ {\cal
Y}_t\right],
\end{equation}
where
$\sigma$-field
${\overline{\YY}}_{t,t-1}=\sigma(\{(Y_s, Y^{2}_{r}) ,\, 1\leq s\leq t,\, 1\leq r\leq t-1\})$.
 Again, conditionally  Gaussian
properties of the pair  $(X,\xi)$ imply that
\begin{equation*}
\frac{\Eg \left[\mathbb{I}(X_{t}  \in dx)\exp\left\{-\xi_{t} \right\}/
{\cal \overline{Y}}_{t,t-1} \right]}{\Eg[\exp(-\xi_{t-1})/
{\cal \overline{Y}}_{t,t-1}]}=
\displaystyle{[2\pi\tg_{t}]^{-\frac{1}{2}}}
\end{equation*}
\begin{equation}\label{gausden}
\times \exp \displaystyle{
\left\{-\frac{1}{2}(x-\widetilde{Z}^{h}_{t})^{2}\tg^{-1}_t
\right\}\,dx},
\end{equation}
where $\displaystyle{\widetilde{Z}^{h}_{t}=\Eg[X_t/
{\bar{\YY}}_{t,t-1}]- \Eg[(X_t-\Eg[X_t/
{\bar{\YY}}_{t,t-1}])(\xi_{t-1}-\bp_{t-1}(\xi_{t-1}))/
{\bar{\YY}}_{t,t-1}]}$ and $\displaystyle{\tg_{t}=\Eg[(X_t-\Eg[X_t/
{\bar{\YY}}_{t,t-1}])]^{2}}$.
Now the  desired equality \eqref{bayesford} follows directly from Proposition~\ref{p:CM}.\\
It is worth emphasizing that (for negative $\mu$)  now we know the probabilistic interpretation of the involved processes
 $(Z^h,\widetilde{Z}^{h}, \,\bar{\gamma},\tg)$. Actually we have proved that $Z^{h}$
is the difference $\bar{\pi}_{t}(X)-\bar{\gamma}_{_{X\xi}}(t)$ and
$\bar{\gamma}$ is nothing but the covariance of the filtering
error of $X$ in view of auxiliary observations $\bar{Y}$.\\ For the pair $(\widetilde{Z}^{h}, \,\tg)$ we have the same relations but with respect to the $\sigma$-field
${\cal
Y}_{t,t-1}=\sigma(\{(Y_s, Y^{2}_{r}) ,\, 1\leq s\leq t,\, 1\leq r\leq t-1\})$.
Of course, after a simple integration of  $\lambda_{t}$, formula (\ref{cdensg})  gives
Proposition~\ref{CMC}  and therefore  the solution of the LEG and RS filtering problems.
Let us also observe that the relations 
$\displaystyle{\tilde{Z}^{h}_{t}=\frac{Z_{t}^h+A_t \bg_t Y_t}{1+ A_t^2 \bg_t }},\,
\displaystyle{\tilde{\gamma}_{t}=\frac{ \bg_t}{1+ A_t^2 \bg_t }}$ which were announced in Remark~\ref{probinterp'}  follow from  the Note following
  Theorem 13.1 in \cite{lipshi1bis}.

\medskip
\textit{Linear Exponential Gaussian Control}

\medskip\noindent

In the context of the LEG  control problem for a partially observed process, the information state is also defined (see, \textit{e.g.}, \cite{elliott1}) as the density $\lambda_{t}$, with respect to the Lebesgue measure,
of the non normalized random measure
$\omega_{t}$:
\begin{equation}\label{mescondcon}
 \omega_{t}(dx)=\displaystyle{\Eg\left[\left.\mathbb{I}(X_t \in dx)\exp\left\{\frac{\mu}{2} \sum_{s=1}^{t-1}
X_{s}^{2}Q_{s} \right\} \right/ {\cal
Y}_t\right]},
\end{equation}
where  $X$ is the controlled state governed  by the equation:
\begin{equation}\label{SDECONT}
X_t=a_tX_{t-1}+b_tu_t +\widetilde{\varepsilon}_t\,,\;t\geq 1\,;\; X_0=0\,,
\end{equation}
 $( \widetilde{\varepsilon}_t)_{t \ge 1}$ is a sequence of i.i.d. standard
Gaussian variables
and $u_{t}\in \YY_{t-1} $ corresponding to the available observation $Y$  defined by the equation \eqref{observ}.\\
By the same way that we have just explained, for  the conditionally Gaussian pair $(X,Y)$, one can check that the density $\lambda_{t}$
 satisfies the following equality:
\begin{multline}\label{cdensgcon}
\lambda_{t}(x) =
\displaystyle{\frac{1}{\sqrt{2\pi \tg_t}}
\exp\left\{-\frac{(x-\widetilde{Z}_{t})^{2}}{2\widetilde{\gamma}_{t}}
 \right\}}\times
\\
\prod_{r=1}^{t-1} \left[\frac{1+S_r \bg_r}{1+A_r^{2} \bg_r}\right]^{-1/2} \times \exp\left\{\frac{\mu}{2} Q_r \frac{1+A_r^2 \bg_r}{1+ S_r \bg_r} \times \widetilde{Z}_{t}^{2}\right\} \times \mathcal{M}_t,
\end{multline}
where $\displaystyle{\tilde{\gamma}_{t}=\frac{ \bg_t}{1+ A_t^2 \bg_t }}$, $\bar{\gamma}$ is the solutions of
equation (\ref{GAMMABAR}), the martingale  $(\mathcal{M}_t)_{\ge 1}$ is defined by \eqref{martin}
 and $\widetilde{Z}$ is the solution of the equation
\begin{equation}\label{z:repres}
\widetilde{Z}_{t} = \frac{a_t}{1+ S_t \bg_t }\widetilde{Z}_{t-1}+\frac{b_t}{1+ S_t \bg_t }u_t
  + \bar{\gamma}_t A_tY_{t}.
\end{equation}
Actually it is the equation for the difference $\widetilde{Z}=\bar{\pi}_{t,t-1}(X)-\bar{\gamma}_{X\xi}(t,t-1)$, where the conditional expectations are taken with respect to the auxiliary observation process $\bar{Y}$ defined by the equations \eqref{observ} and \eqref{Yaux} with $h=0$.\\
Equality \eqref{cdensgcon} gives the possibility to rewrite the cost function in terms of the completely  observable process $\widetilde{Z}$, namely:
\begin{equation}\nonumber
\begin{array}{ccl}
\Eg\Big[
\exp\Big\{\displaystyle{\frac{\mu}{2} \sum_{s=1}^T
X_{s}^{2}Q_{s} }
\Big]
=\Eg\Big\{ \Eg\Big[\left.
\exp\Big\{\displaystyle{\frac{\mu}{2} \sum_{s=1}^T
X_{s}^{2}Q_{s} }\Big\} \right/
{\cal Y}_T\Big]\Big\}
\\
=\displaystyle{\prod_{r=1}^{t-1} \left[\frac{1+S_r \bg_r}{1+A_r^{2} \bg_r}\right]^{-1/2}} \Eg\Big[\displaystyle{\exp\left\{\frac{\mu}{2}\sum_{s=1}^T\widetilde{Z}_{s}^{2}\widetilde{Q}_{s}
\right\}}\Big]\times \mathcal{M}_T
\\
=\displaystyle{\prod_{r=1}^{t-1} \left[\frac{1+S_r \bg_r}{1+A_r^{2} \bg_r}\right]^{-1/2}} \widetilde{\Eg}\Big[\displaystyle{\exp\left\{\frac{\mu}{2}\sum_{s=1}^T\widetilde{Z}_{s}^{2}\widetilde{Q}_{s}
\right\}}\Big],
\end{array}
\end{equation}
where $\displaystyle{\widetilde{Q}_{r} = Q_r \frac{1+A_r^2 \bg_r}{1+ S_r \bg_r}}$ and $\widetilde{\Eg}$ stands for an expectation with respect to the new measure $\widetilde{\Pg}$ such that:
$$
\frac{d\widetilde{\Pg}}{d\Pg} =\mathcal{M}_T .
$$
With respect to this new measure the solution of  equation \eqref{z:repres} can be represented as
\begin{equation}\label{z:repres'}
\widetilde{Z}_{t} = a_t \frac{1+A_t^2 \bg_t}{1+ S_t \bg_t}\widetilde{Z}_{t-1}+b_t \frac{1+A_t^2 \bg_t}{1+ S_t \bg_t}u_t
  + \frac{\bar{\gamma}_tA_t}{1+A_r^{2} \bg_r}\bar{\varepsilon}_{t},
\end{equation}
where $( \bar{\varepsilon}_t)_{t \ge 1}$ is a new sequence of i.i.d. standard
Gaussian variables.
Thus, the new  process $\widetilde{Z}$ plays the role of the  completely observed controlled state (see \cite{bensoussan} and \cite{elliott1}).

Now we emphasize that the probabilistic interpretation of the
``information state" $\widetilde{Z}$, used in \cite{elliott1}
is nothing but
$\widetilde{Z}_{t}=\bar{\pi}_{t,t-1}(X)-\bar{\gamma}_{X\xi}(t)$, where the conditional expectations are taken with respect to the auxiliary observation process $\bar{Y}$ defined by the equations \eqref{observ} and \eqref{Yaux} with $h=0$. Also, $\bar{\gamma}$ is the conditional covariance of $X$.

\section{Complementary part - More general case}\label{complements}

In this section we analyze LEG and RS filtering problems in a more general contexts when we do not suppose a special structure of the observation sequence $(Y_t)_{t\ge 1}$. We suppose only that the process $(X_t,\,Y_t)_{t\ge 1}$ is Gaussian (even  conditionally Gaussian). Our goal is to reduce LEG (RS) filtering problems to an auxiliary  risk-neutral filtering problem.
First of all we fix $\mu=-1$ and we will find the probabilistic interpretation of the solution. After to find the solution for $\mu \ne -1$ we shall have only to replace $Q$ by $-\mu Q$ in the answer.
So, let $(Y^{2}_{t},\, \xi_{t})$ be defined by equations \eqref{Yaux} - \eqref{xi:eq}  and let us denote by
\begin{equation}\label{Ztildainterpr}
\widetilde{Z}^{h}_{t}=\Eg[X_t/
{\bar{\YY}}_{t,t-1}]- \Eg[(X_t-\Eg[X_t/
{\bar{\YY}}_{t,t-1}])(\xi_{t-1}-\bp_{t-1}(\xi_{t-1}))/
{\bar{\YY}}_{t,t-1}],
\end{equation}
\begin{equation}
\tg_{t}=\Eg[(X_t-\Eg[X_t/
{\bar{\YY}}_{t,t-1}])]^{2},
\end{equation}
where ${
\bar{\YY}}_{t,t-1}$ is the
$\sigma$-field
${
\bar{\YY}}_{t,t-1}=\sigma(\{(Y_s, Y^{2}_{r}) ,\, 1\leq s\leq t,\, 1\leq r\leq t-1\})$.
Again, let $\displaystyle {J_t= \exp\left\{-\frac{1}{2} \sum_{s=1}^t (X_s-h_s)^2 Q_s \right\}}$ and let us denote by $\mI_t$ the conditional expectation
$\displaystyle {\mI_t= \pi_t(J_t)}$, or
$$
\mI_t= \E \left(\left.\exp\left\{-\frac{1}{2} \sum_{s=1}^t (X_s-h_s)^2 Q_s \right\} \right/ {\YY}_{t}\right),
$$
where $h_s\in \YY_{s}, \, s\ge 1$.
We claim the following generalization of Proposition~\ref{p:CM}.

\begin{proposition}\label{p:CMG}
The following equality holds for any $T\ge 1$:
$$
\mI_T=\prod_{t=1}^T \left[1+Q_t\tg_{t}\right]^{1/2} \times \exp\left\{-\frac{1}{2} \frac{Q_t}{1+ Q_t \tg_t} \times \left[ h_t -\widetilde{Z}^{h}_{t}  \right]^{2}\right\} \times \mathcal{M}_T,
$$
where   $(\mathcal{M}_T)_{T\ge 1}$ is a martingale defined by :
\begin{multline}\label{martin}
\mathcal{M}_T = \prod_{t=1}^T \left[\frac{\sigma_t^{2}}{\bar{\sigma}_t^{2} }\right]^{1/2} \exp\left\{ \frac{1}{2\sigma_t^{2}} \, (Y_t - \pi_{t-1}(Y_t) )^{2} -\frac{1}{2\bar{\sigma}_t^{2}} \, (Y_t -\bar{ V}_{t} )^{2}
\right\},
\end{multline}
where
$$
\sigma_t^{2}=\Eg(Y_t-\pi_{t-1}(Y_t))^{2},\,
\bar{\sigma}_t^{2}=\Eg(Y_t-\bar{\pi}_{t-1}(Y_t))^{2},
$$
$$
\bar{ V}_{t}=\bar{\pi}_{t-1}(Y_t) -\bg_{_{Y\xi}}(t),\, \bg_{_{Y\xi}}(t)=\Eg[(Y_t-\bar{\pi}_{t-1}(Y_t))
(\xi_{t-1}-\bp_{t-1}(\xi_{t-1}))/
{\bar{\YY}}_{t-1}].
$$
\end{proposition}

\paragraph{Proof}To prove Proposition~\ref{p:CMG}
 let us again calculate the ratio

$$
\frac{\mI_t}{\mI_{t-1}} = \frac{\pi_t(J_t)}{\pi_{t-1}(J_{t-1})}= \frac{\pi_t(J_t)}{\pi_{t}(J_{t-1})}\frac{\pi_{t}(J_{t-1})}{\pi_{t-1}(J_{t-1})}=
$$
$$
=\frac{\pi_t(J_t)}{\pi_{t}(J_{t-1})}\frac{\mathcal{M}_t}{\mathcal{M}_{t-1}}
$$
with a martingale $\mathcal{M}_t$ such that:
\begin{equation}\label{martgen}
\mathcal{M}_{t}=\prod_{s=1}^{t}\frac{\pi_{s}(J_{s-1})}{\pi_{s-1}(J_{s-1})}.
\end{equation}
The same arguments that we used in the proof of Proposition~\ref{p:CM} show that
$$
\frac{\pi_t(J_t)}{\pi_{t}(J_{t-1})}= \frac{\bp_{t,t-1}(\exp \{-\frac{1}{2}Q_t (X_t-h_t)^2 - \xi_{t-1}\} )}{\bp_{t,t-1}(\exp(-\xi_{t-1}))}
$$
$$
=(1+Q_t\tg(t))^{-1/2}  \exp\left\{ - \frac{1}{2} \cdot\frac{Q_t}{1+Q_t\tg(t)}(\widetilde{Z}^{h}_{t}-h_t)^2 \right\}.
$$
To finish the proof  we turn to the representation of the martingale $\mathcal{M}_t$ defined by \eqref{martgen}.
First of all we claim that
\begin{equation}\label{martbayes}
\frac{\mathcal{M}_{t}}{\mathcal{M}_{t-1}}=\left.\frac{\widetilde{\pi}_{t-1}(\mathbb{I}(Y_t \in dy))}{\pi_{t-1}(\mathbb{I}(Y_t \in dy))}\right|_{y=Y_t},
\end{equation}
where $\widetilde{\pi}$ stands for the conditional expectation with respect to the measure $\widetilde{\Pg}$ such that
$\displaystyle{\frac{d \widetilde{\Pg}}{d \Pg}={\mathcal{M}_{T}}}$. Indeed, it is the direct consequence of the classical Bayes formula
$$
\widetilde{\pi}_{t-1}(\mathbb{I}(Y_t \in dy))=\frac{\pi_{t-1}(\mathbb{I}(Y_t \in dy){\mathcal{M}_{T}})}{\pi_{t-1}({\mathcal{M}_{t-1}})}
=\pi_{t-1}(\mathbb{I}(Y_t \in dy){\mathcal{M}_{t}}).
$$
To finish the proof it is sufficient to note that representations \eqref{martgen} and \eqref{martbayes} imply that
$$
\frac{\mathcal{M}_t}{\mathcal{M}_{t-1}}=\left.\frac{\pi_{t}(J_{t-1})}{\pi_{t-1}(J_{t-1})}=\left.\frac{\widetilde{\pi}_{t-1}(\mathbb{I}(Y_t \in dy))}{\pi_{t-1}(\mathbb{I}(Y_t \in dy))}\right|_{y=Y_t}=\frac{\pi_{t-1}(\mathbb{I}(Y_t \in dy)\pi_{t}(J_{t-1}))}{\pi_{t-1}(J_{t-1})\pi_{t-1}(\mathbb{I}(Y_t \in dy))}\right|_{y=Y_t}=
$$
$$
\left.\frac{\pi_{t-1}(\mathbb{I}(Y_t \in dy)J_{t-1})}{\pi_{t-1}(J_{t-1})\pi_{t-1}(\mathbb{I}(Y_t \in dy))}\right|_{y=Y_t}.
$$
Again, we can use the same arguments that we used in the proof of Proposition~\ref{p:CM}:
$$
\frac{\pi_{t-1}(\mathbb{I}(Y_t \in dy)J_{t-1})}{\pi_{t-1}(J_{t-1})}=\frac{\bp_{t-1}(\mathbb{I}(Y_t \in dy)\exp\{-\xi_{t-1}\})}{\bp_{t-1}(\exp\{-\xi_{t-1}\})}=
$$
$$
=\frac{1}{\sqrt{2\pi \bar{\sigma}^{2}}}\exp\left(-\frac{(Y_t-\bar{V}_t )^{2}}{2\bar{\sigma}^{2}}\right).
$$

A direct consequence of Proposition~\ref{p:CMG} is the following statement:

\begin{corollary}
Let $\bh$ be the solution of LEG (and RS) filtering problem \eqref{LEGdef} (and \eqref{rde}). Then the following equality holds for any $t \ge 1$:

$$
\bh_t= \widetilde{Z}^{\bh}_{t}.
$$
\end{corollary}

\section{Particular cases - again}\label{againpart}
\subsection{Markov type observations}
Here we turn to the case when the  observations $(Y_t)_{t\ge 1}$ are conditionally independent given $X$. More precisely,
 we deal with a signal-observation model
$
(X_t,Y_t)_{t\ge 1},
$
where the signal $X=(X_t)_{t\ge 1},\, X_t \in \mathbb{R}^{n}$ is an arbitrary Gaussian sequence with mean vector
$m=(m_t, t\geq 1)$ and covariance matrix  $K =(K(t,s), t\geq 1,
s\geq 1)$, \textit{i.e.},
$$
\Eg X_t=m_t,\quad\Eg (X_t-m_t)(X_s-m_s)^{\prime}=K(t,s)\,.
t\geq 1\,,\; s\geq 1\,,
$$
The observation process $Y=(Y_t,\, t\ge 1)$ is given by
\begin{equation}\label{obser}
Y_t= A_t X_t +\varepsilon_t,
\end{equation}
for some sequence $A=(A_t,\, t\ge 1)$ of $m\times n$ matrices, where $\varepsilon=(\varepsilon_{t})_{t\ge 1}$ is a sequence of i.i.d. $\mN(0,Id)$ random variables and $\varepsilon$ and $X$ are independent.
In this case  we can write the multidimensional analogue  of the equation \eqref{zhtild}, which is nothing else but the dynamic equation for the process $\widetilde{Z}^{h}$ defined by \eqref{Ztildainterpr}. We obtain:
$$
 \tilde{Z}^{h}_{t}=m_t +\sum_{l=1}^{t-1}\bg(t,l)[Id+\bg_{l}(A^{\prime}_{l}A_{l}-\mu Q_{l})]^{-1}\mu Q_{l} (h_l - \tilde{Z}_{l}^h) + \sum_{l=1}^{t}\bg(t,l) A^{\prime}_{l} (Y_l - A_l \tilde{Z}_{l}^{h}),
$$
where the matrix $\bg(t,l)$ satisfies the following equation (which is the multidimensional analog of the equation \eqref{GAMMABAR}):
\begin{equation}\label{gammamultdim}
\bg(t,s)=K(t,s)-\sum_{l=1}^{s-1} \bg(t,l)\bar{A}_{l}^{\prime}[Id +\bar{A}_{l}\bg_l\bar{A}_{l}^{\prime}]^{-1}\bar{A}_{l}^{\prime} \bg^{\prime}(s,l),
\end{equation}
where $\bar{A}_{l}=\left(%
\begin{array}{c}
  A_l \\
  -\mu Q_l \\
\end{array}%
\right).$

Now the solution of the LEG (and RS) filtering problem $\bar{h}$ is nothing else but:
$$
\bar{h}_{t}=m_t + \sum_{l=1}^{t}\bg(t,l) A^{\prime}_{l} (Y_l - A_l\bar{h}_{l}).
$$
\subsection{Markov type observations, correlated signal and observation noises}
Let us drop the assumption  that $X$ and $\varepsilon $ in the observation equation \eqref{obser} are independent. Denote by $K_{_{X\varepsilon}}(t,s)$ the covariance matrix of the signal and the observation noise,  \textit{i.e.},
$$
\Eg (X_t-m_t)\varepsilon_s^{\prime}=K_{_{X\varepsilon}}(t,s), \quad
t\geq 1\,,\; s\geq 1.
$$
It can be checked that the following slight modification of the previous statement holds. \\ Let the matrix $\bg(t,l)$ be the unique solution  of the following equation
\begin{multline}\label{gammacorrel}
\bg(t,s)=K(t,s)-\sum_{l=1}^{s-1} [\bg(t,l)\bar{A}_{l}^{\prime}+\bar{K}_{_{X\varepsilon}}(t,l)]
\\
[Id +\bar{A}_{l}\bg_l\bar{A}_{l}^{\prime} +\bar{A}_{l}\bar{K}_{_{X\varepsilon}}(l,l)+\bar{K}_{_{X\varepsilon}}(l,l)^{\prime}\bar{A}_{l}^{\prime}]^{-1}
\\
[\bar{A}_{l}^{\prime} \bg^{\prime}(s,l)+\bar{K^{\prime}}_{_{X\varepsilon}}(s,l)],
\end{multline}
with $\bar{A}_{l}=\left(%
\begin{array}{c}
  A_l \\
  -\mu Q_l \\
\end{array}%
\right),\quad \bar{K}_{_{X\varepsilon}}(t,l)=(K_{_{X\varepsilon}}(t,l)\quad  \mathbf{0}).$\\
Then the solution of the LEG (and RS) filtering problem $\bar{h}$ satisfies the following equation
\begin{equation}\label{hbarcorrel}
\bar{h}_{t}=m_t + \sum_{l=1}^{t}[\bar{K}_{_{X\varepsilon}}(t,l)+\bg(t,l) A^{\prime}_{l}][Id +A_{l}K_{_{X\varepsilon}}(l,l)]^{-1} (Y_l - A_l\bar{h}_{l}).
\end{equation}

\subsection{Observations containing Moving Averages of order~1}\label{MA1XY}
Now we consider the case of a  MA(1) type process, {\em i.e.}, the following signal-observation model:
$$
X_{t} = \widetilde{\varepsilon}_{t} + \lambda \widetilde{\varepsilon}_{t-1}\,; t\ge 1\,,
$$
$$
Y_{t} =\alpha_t X_t + \varepsilon_{t} + \beta \varepsilon_{t-1}\,; t\ge 1\,,
$$where $( \varepsilon_t, \widetilde{\varepsilon}_t)_{t\ge 0}$~is a sequence of i.i.d.
Gaussian variables and $\lambda$  and $\beta$ are  real numbers.\\
Let us denote by $A_t$ the row $\bar{A}_{t}=(\alpha_{t}\quad \beta)$ and by $\bar{X}_{t}$  the vector  $\bar{X}_{t}=\left(%
\begin{array}{c}
  X_t \\
  \varepsilon_{t-1} \\
\end{array}%
\right).$
 Of course
$\bar{X}$ is centered, has the
covariance matrix
$K(t,s) =\left(%
\begin{array}{cc}
 (1+\lambda^{2}) \mathbf{1}(s=t-1)+ \lambda \mathbf{1}(s=t)& 0 \\
  0 & \mathbf{1}(s=t) \\
\end{array}%
\right)$
and the covariance between $\bar{X}$ and $\varepsilon$ is $K_{_{X\varepsilon}}(t,s)=\left(%
\begin{array}{c}
  0 \\
 \mathbf{1}(s=t-1)  \\
\end{array}%
\right).$

The solution $\bg$ to  (\ref{gammacorrel}) then can be found as:
$$
\bg(t,s) =\mathbf{0}\,,\; s < t-1\,; \quad  \bg(t,t-1) = \left(%
\begin{array}{cc}
  \lambda & 0 \\
  0 & 0 \\
\end{array}%
\right)\,,\;
t\ge 1\,,
$$
and $\bg(t,t)=\bg_{t}$ where $\bg_{t}$  is
the solution of the equation:
\begin{multline*}
\bg_{t} =\left(%
\begin{array}{cc}
  1+\lambda^{2} & 0 \\
  0 & 1 \\
\end{array}%
\right) +
\\
 +\left(%
\begin{array}{cc}
  \lambda\alpha_{t-1} & -\lambda\mu Q_{t-1}\\
  0 & 1 \\
\end{array}%
\right)\left[ Id + \left(%
\begin{array}{cc}
  \alpha_{t-1} & \beta \\
  -\mu Q_{t-1} & 0 \\
\end{array}%
\right)\bg_{t-1}\left(%
\begin{array}{cc}
  \alpha_{t-1} &  -\mu Q_{t-1} \\
  \beta & 0 \\
\end{array}%
\right)\right]^{-1}
\\
\times \left(%
\begin{array}{cc}
  \lambda\alpha_{t-1} & 0 \\
  -\lambda\mu Q_{t-1} & 1 \\
\end{array}%
\right) \,,\;t\ge 1\,;
 \quad \bg_0= \mathbf{0}.
\end{multline*}
provided that this equation has the unique nonnegative definite solution.

Moreover, equation
\eqref{hbarcorrel} for the solution $\bh$ of the LEG filtering problem (\ref{LEGdef}) can be reduced to the
following one:
\begin{equation*}
\bh_t=\Lambda_{t}^{-1} \left(%
\begin{array}{c}
 \lambda\alpha_{t}  \\
  1+\beta \\
\end{array}%
\right) [Y_{t-1}-A_{t-1}\bh_{t-1}] +\Lambda_{t}^{-1}\bg_t \left(%
\begin{array}{c}
 \alpha_{t}  \\
  \beta \\
\end{array}%
\right)Y_t, \, t\ge 1,\, \bh_0=0,
\end{equation*}
with $\Lambda_{t}= Id + \bg_{t}A_{t}^{\prime}A_{t}.$

\subsection{Observations containing Gaussian AR(1) process }\label{GMSXY}
In this part we concentrate on the case of a
Gaussian AR(1) type process $Y$, {\em i.e.},
\begin{equation}\label{model ARY}
Y_{t}= \alpha_{t} X_{t}+\varepsilon_{t}\,,\; t\ge 1\,;
\quad Y_{0}=0\,,
\end{equation}
where
$$
\varepsilon_{t}=b \varepsilon_{t-1}+\widetilde\varepsilon_{t},
$$
and  $(\widetilde\varepsilon_{t}, \; t=1,2,\dots)$ is a sequence of i.i.d.
standard Gaussian random variables independent of $X$.
We also suppose that the signal $X$ is a Gaussian AR(1) process,
{\em i.e.},
\begin{equation}\label{model ARY}
Y_{t}= a_{t} X_{t}+\epsilon_{t}\,,\; t\ge 1\,;
\quad X_{0}=0\,,
\end{equation}
and also $(\epsilon_{t}, \; t=1,2,\dots)$ is a sequence of i.i.d.
standard Gaussian random variables.
Proceeding  as in Sections \ref{GMS} and \ref{MA1XY} we can write the dynamic equation for the solution of LEG and RS filtering problems $\bar{h}$. Namely, $\bar{h}$ is the first component $\bar{h}^{1} $ of the solution of the following recursive equation:
\begin{equation*}
\bh_t=\bg_t A_t^{\prime}[Y_{t}-A_{t}\bh_{t}]  +\left(%
\begin{array}{cc}
a_{t} & 0 \\
 0 & b \\
\end{array}%
\right) \bh_{t-1} +  \left(%
\begin{array}{c}
 0  \\
  1 \\
\end{array}%
\right)[Y_{t-1}-A_{t-1}\bh_{t-1}] \, t\ge 1,\, \bh_0=\mathbf{0},
\end{equation*}
where  $A_t= (\alpha_{t} \quad \beta)$ and $\bg$ is the unique nonnegative defined solution of the Ricatti equation:
\begin{multline*}
\bg_t=\left(%
\begin{array}{cc}
  a_t & 0 \\
  0 & b\\
\end{array}%
\right)\bg_{t-1}\left(%
\begin{array}{cc}
  a_t & 0 \\
  0 & b\\
\end{array}%
\right) -\left[\left(%
\begin{array}{cc}
  a_t & 0 \\
  0 & b\\
\end{array}%
\right)\bg_{t-1}\bar{A}^{\prime}_{s-1} +\left(%
\begin{array}{ccc}
  0 & 0 & 0\\
  0 & 1 & 0\\
\end{array}%
\right)\right]
\\
\times\left[Id +\bar{A}_{t-1} \bg_{t-1}\bar{A}^{\prime}_{s-1}\right]^{-1}\left[\bar{A}_{t-1}\bg_{t-1}\left(%
\begin{array}{cc}
  a_t & 0 \\
  0 & b\\
\end{array}%
\right) +\left(%
\begin{array}{cc}
  0 & 0 \\
  0 & 1 \\
  0 & 0 \\
\end{array}%
\right)\right],
\end{multline*}
with $\bar{A}_{t}=\left(%
\begin{array}{cc}
  \alpha_{t} & b \\
  -\mu Q_t & 0 \\
  0 & 0 \\
\end{array}%
\right).$
\bibliographystyle{plain}

\end{document}